\documentclass[12pt]{article}

\usepackage{geometry} 

\newtheorem{tm}{Theorem}

\newtheorem{kor}{Corollary}
\newtheorem{rem}{Remark}

\usepackage{amssymb}
\usepackage{amsmath,amsfonts,latexsym}
\usepackage{latexsym}
\usepackage{float}

\begin{document}

\vspace*{0.5cm}

\begin{center}
{\Large\bf  On transitive designs and strongly regular graphs constructed from Mathieu group $M_{11}$}
\end{center}

\vspace*{0.5cm}

\begin{center}
        Dean Crnkovi\'c (deanc@math.uniri.hr)\\
		and \\
		Andrea \v Svob (asvob@math.uniri.hr)\\[3pt]
		{\it\small Department of Mathematics} \\
		{\it\small University of Rijeka} \\
		{\it\small Radmile Matej\v ci\'c 2, 51000 Rijeka, Croatia}\\
\end{center}

\vspace*{0.5cm}

\begin{abstract}
In this paper we construct structures from Mathieu group $M_{11}$. We classify transitive $t$-designs with 11, 12 and 22 points admitting a transitive action of Mathieu group $M_{11}$. Thereby we proved the existence of designs with parameters $3$-$(22,7,18)$ and found first simple designs with parameters $4$-$(11,5,6)$ and $5$-$(12,6,6)$. Additionally, we proved the existence of $2$-designs with certain parameters having 55 and 66 points. Furthermore, we classified strongly regular graphs on at most 450 vertices admitting a transitive action of the Mathieu group $M_{11}$.
\end{abstract}

\bigskip

{\bf AMS classification numbers:} 05B05, 20D06, 05E18, 05E30.


\section{Introduction}

We assume that the reader is familiar with the basic facts of group theory, design theory and theory of strongly regular graphs.
We refer the reader to \cite{bjl, tonchev-book} for relevant background reading in design theory, to \cite{atlas, r} for relevant background reading in group theory, and to 
background reading in theory of strongly regular graphs we refer the reader to \cite{bjl, crc-srg, tonchev-book}.

\bigskip

An incidence structure is an ordered triple
${\mathcal{D}}=({\mathcal{P}},{\mathcal{B}},{\mathcal{I}})$ where
${\mathcal{P}}$ and ${\mathcal{B}}$ are non-empty disjoint sets
and $ {\mathcal{I}}\subseteq {\mathcal{P}}\times {\mathcal{B}}$.
The elements of the set ${\mathcal{P}}$ are called points, the
elements of the set ${\mathcal{B}}$ are called blocks and
$\mathcal{I}$ is called an incidence relation.
If $|{\mathcal{P}}|=|{\mathcal{B}}|$, then the incidence structure is
called symmetric. The incidence
matrix of an incidence structure is a $v \times b$ matrix
$[m_{ij}]$ where $v$ and $b$ are the numbers of points and blocks
respectively, such that $m_{ij} = 1$ if the point $P_i$ and the block 
$x_j$ are incident, and $m_{ij} = 0$ otherwise. An isomorphism
from one incidence structure to another is a bijective mapping of
points to points and blocks to blocks which preserves incidence.
An isomorphism from an incidence structure ${\mathcal{D}}$ onto
itself is called an automorphism of ${\mathcal{D}}$. The set of
all automorphisms forms a group called the full automorphism group
of ${\mathcal{D}}$ and is denoted by $Aut({\mathcal{D}})$.\\

A $t$-$(v,k,\lambda)$ design is a finite incidence structure
${\mathcal{D}}=({\mathcal{P}},{\mathcal{B}},{\mathcal{I}})$ satisfying the
following requirements:
\begin{enumerate}
  \item $|{\mathcal{P}}|=v$,
  \item every element of ${\mathcal{B}}$ is incident with exactly
  $k$ elements of ${\mathcal{P}}$,
  \item every $t$ elements of ${\mathcal{P}}$ are incident with exactly
$\lambda$ elements of ${\mathcal{B}}$.
\end{enumerate}
If ${\mathcal{D}}$ is a $t$-design, then it is also an $s$-design, for $1 \leq s \leq t-1$. A $2$-$(v,k, \lambda)$ design is called a block design. 
We say that a $t$-$(v,k,\lambda)$ design $\mathcal{D}$ is a quasi-symmetric design with intersection numbers $x$ and $y$ $(x<y)$ 
if any two blocks of $\mathcal{D}$ intersect in either $x$ or $y$ points.

A graph is regular if all the vertices have the same degree; a regular graph is strongly regular of type $(v,k, \lambda , \mu )$ if it has $v$ vertices, degree $k$,
and if any two adjacent vertices are together adjacent to $\lambda$ vertices, while any two non-adjacent vertices are together adjacent to $\mu$ vertices.
A strongly regular graph  of type $(v,k, \lambda , \mu )$ is usually denoted by SRG$(v, k, \lambda, \mu)$.

\bigskip

In this paper we consider $t$-designs and strongly regular
graphs constructed from the Mathieu group $M_{11}$. 
It is the simple group of order $7920$, and up to conjugation it has $39$
subgroups given in Table \ref{tb:subgrpsM11}.

Using the method introduced in \cite{cms}, we classify all $t$-designs  on $11$, $12$ or $22$ points on which the group $M_{11}$ 
acts transitively on points and blocks.
Additionally, we obtained numerous transitive designs, under the action of $M_{11}$, for $v=55, 66$. In many cases we proved the existence of $2$-designs with certain parameters. We also proved the existence of $4$-$(11,5,6)$, $5$-$(12,6,6)$ and $3$-$(22,7,18)$ designs.

Further, we construct strongly regular graphs on 55, 66, 144 or 330 vertices from the simple group $M_{11}$. Constructed strongly regular graphs have been known before.

Generators of the group $M_{11}$ are available on the Internet:\\ http://brauer.maths.qmul.ac.uk/Atlas/. 
All the structures are obtained by using programmes written for Magma \cite{magma}.\\

The paper is organized as follows: in Section \ref{constr} we
briefly describe the method of construction of transitive designs used in this paper, and in Section \ref{results} we describe $t$-designs constructed under the action of the Mathieu group $M_{11}$.

\section{Structures constructed from groups} \label{constr}

The  construction of primitive symmetric $1$-designs and regular graphs for which the stabilizer 
of a point and the stabilizer of a block are conjugate
is presented in \cite{km}, \cite{km1} and \cite{key1}. The generalization, 
{\it i.e.} the method for constructing not necessarily symmetric but still primitive $1$-designs, 
is presented in \cite{cm1} and \cite{cm}. 
In \cite{cms}, a construction of not necessarily primitive, but still transitive block designs is presented.

\begin{tm}[\cite{cms}]\label{main}
Let $G$ be a finite permutation group acting transitively on the sets $\Omega_1$ and $\Omega_2$
of size $m$ and $n$, respectively.
Let $\alpha \in \Omega_1$ and $\Delta_2 =  \bigcup_{i=1}^s \delta_i G_{\alpha}$, where $G_{\alpha} = \{ g \in G \ | \ \alpha g = \alpha \}$ is the stabilizer of $\alpha$
and $\delta_1,...,\delta_s \in \Omega_2$ are representatives of distinct $G_\alpha$-orbits on $\Omega_2$.
If $\Delta_2 \neq \Omega_2$ and
$${\mathcal{B}}=\{ \Delta_2 g : g \in G \},$$
then ${\mathcal{D}}(G,\alpha,\delta_1,...,\delta_s)=(\Omega_2,{\mathcal{B}})$ is a
$1$-$(n, | \Delta_2 |, \frac{|G_{\alpha}|}{|G_{\Delta_2}|}\sum_{i=1}^{s} | \alpha G_{\delta_i} |)$ design with $\frac{m\cdot |G_{\alpha}|}{|G_{\Delta_2}|}$ blocks.
The group $H\cong G/{\bigcap_{x\in \Omega_2}G_x}$ acts as an automorphism group on $(\Omega_2,{\mathcal{B}})$, 
transitively on points and blocks of the design.

If $\Delta_2=\Omega_2$ then the set $\mathcal{B}$ consists of one block, and  ${\mathcal{D}}(G,\alpha,\delta_1,...,\delta_s)$ is
a design with parameters $1$-$(n,n,1)$.
\end{tm}

The construction described in Theorem \ref{main} gives us all simple designs on which the group $G$ acts transitively on the points and blocks, {\it i.e.} if $G$ acts 
transitively on the points and blocks of a simple $1$-design ${\mathcal{D}}$, then ${\mathcal{D}}$ can be obtained as described in Theorem \ref{main}. 

If a group $G$ acts transitively on $\Omega$, $\alpha \in \Omega$, and $\Delta$ is
an orbit of $G_{\alpha}$, then 
$\Delta' = \{ \alpha g \ | \ g \in G,\ \alpha {g^{-1}} \in \Delta \}$
is also an orbit of $G_{\alpha}$. $\Delta'$ is called the orbit of $G_{\alpha}$
paired with $\Delta$. It is obvious that $\Delta'' = \Delta$ and $| \Delta' | = | \Delta |$.
If $\Delta' = \Delta$, then $\Delta$ is said to be self-paired.

\begin{kor} \label{main-graph}
If $\Omega_1=\Omega_2$ and $\Delta_2$ is a union of self-paired and mutually paired orbits of $G_{\alpha}$, then the design 
${\mathcal{D}}(G,\alpha,\delta_1,...,\delta_s)$ is a symmetric self-dual design and the incidence matrix of that design is the 
adjacency matrix of a $|\Delta_2|-$regular graph. 
\end{kor}

If a group $G$ acts $t$-transitively on the set $\Omega_2$, then the obtained design $(\Omega_2,{\mathcal{B}})$ is a $t$-design (see \cite{cms}).

Using Theorem \ref{main} and Corollary \ref{main-graph} from \cite{cms}, we construct $t$-designs and strongly regular graphs from Mathieu group $M_{11}$.

The method of constructing designs and regular graphs described in Theorem \ref{main} is a generalization of results presented in \cite{cm, km, km1}. 
Using Corollary \ref{main-graph}, one can construct all regular graphs admitting a transitive action of the group $G$, but we will be interested only in those 
regular graphs that are strongly regular.

\section{Structures from $M_{11}$}\label{results}

The Mathieu group $M_{11}$ is a simple group of order 7920, and up to conjugation it has 39 subgroups. It is the smallest sporadic group and acts $4$-transitively on 11 points. There are five simple Mathieu groups, introduced by Emile Mathieu in \cite{mathieu1, mathieu2, mathieu3} and $M_{11}$ is the smallest among all Mathieu groups.
In Table \ref{tb:subgrpsM11} we give the list of all the subgroups, up to conjugation and some of their properties.
Since each transitive action of a group $G$ is permutation isomorphic to an action of $G$ on cosets of its subgroup, the indices of
the subgroups in Table \ref{tb:subgrpsM11} give us degrees of all transitive actions of the group $M_{11}$.

\begin{table}[H]
\begin{center} \begin{scriptsize}
\begin{tabular}{|c|c|c||c|c|c|}
\hline
Structure& Order& Index&Structure& Order& Index\\
\hline\hline
$I$ & 1 & 7920 & $Z_5 : Z_4$ &20 &396\\
$Z_2$ & 2& 3960 & $SL(2,3)$ &24 &330\\
$Z_3$ &3 &2640 & $S_4$ &24 &330\\
$Z_5$ &5 & 1584 & $E_9:Z_4$ &36 &220\\
$E_4$ &4& 1980 &$E_9:Z_4$& 36 &220\\
$Z_4$ &4&1980& $S_3\times S_3$ &36&220\\
$S_3$ &6 &1320 & $GL(2,3)$ &48 &165\\
$S_3$ &6&1320 & $Z_{11}:Z_5$ &55 &144\\
$Z_6$ &6 &1320 & $A_5$& 60 &132\\
$Q_8$ &8&990 &   $A_5$ &60 &132\\                 
$D_8$ &8&990 & $E_9:Q_8$ &72 &110\\
$Z_8$ &8 &990 & $(S_3\times S_3):Z_2$ &72&110 \\
$E_9$ &9 &880 & $E_9:Z_8$ &72 &110\\
$D_{10}$ &10 &792 & $S_5$ &120 &66\\
$Z_{11}$ &11 &720 &$(E_9:Z_8):Z_2$ &144&55\\
$A_4$& 12 &660 & $A_6$ &360&22\\
$D_{12}$ &12 &660 & $PSL(2,11)$ &660&12 \\
$QD_{16}$ &16 &495 & $A_6.Z_2$ &720&11\\
$E_9:Z_2$ &18&440 &$M_{11}$ &7920&1\\
$Z_3\times S_3$ &18&440 &&&\\
\hline\hline
\end{tabular}\end{scriptsize} \caption{Subgroups of the group $M_{11}$}\label{tb:subgrpsM11}
\end{center}
\end{table}

\subsection{$t$-designs with $v \leq 22$}

In this section we give all
$t$-designs with at most 22 points on which the group $M_{11}$ acts transitively. The
designs are obtained from the group $M_{11}$ by using Theorem \ref{main}. In that case, the stabilizers of points are subgroups of $M_{11}$ having the indices 11, 12 and 22. 
The list of all designs obtained is given in Table \ref{tb:td22}. 
In each table we give the parameters of the constructed structures, the number of non-isomorphic structures and their full automorphism group.
The group $M_{11}$ acts 4-transitively on 11 points, hence all designs obtained by Theorem \ref{main} on 11 points are 4-designs.

\begin{table}[H]
\begin{center} \begin{footnotesize}
\begin{tabular}{|c|c|c|c|}
\hline
Parameters of designs & \# of blocks & \# non-isomorphic  & Full automorphism group  \\
\hline
\hline
$3$-$(11, 3, 1)$ & 165 & 1 & $S_{11}$\\ \hline
$4$-$(11, 4, 1)$ & 330 & 1 & $S_{11}$\\ \hline
$4$-$(11, 5, 1)$ & 66 & 1 &  $M_{11}$\\ \hline
$4$-$(11, 5, 6)$ & 396 & 1 & $M_{11}$\\ \hline
\hline
$3$-$(12, 3, 1)$ & 220 & 1 & $S_{12}$\\ \hline
$3$-$(12, 4, 6)$ & 330 & 1 & $M_{11}$\\ \hline
$3$-$(12, 4, 3)$ & 165 & 1 &  $M_{11}$\\ \hline
$3$-$(12, 5, 6)$ & 132 & 1 & $M_{11}$\\ \hline
$3$-$(12, 5, 30)$ & 660 & 1 & $M_{11}$\\ \hline
$3$-$(12, 6, 2)$ & 22 & 1 & $M_{11}$\\ \hline
$3$-$(12, 6, 10)$ & 110 & 1 &  $M_{11}$\\ \hline
$5$-$(12, 6, 6)$ & 792 & 1 & $M_{12}$\\ \hline
\hline
$2$-$(22, 7, 36)$ & 396 & 1 & $M_{11}$\\ \hline
$2$-$(22, 7, 180)$ & 1980 & 1 & $M_{11}$\\ \hline
$2$-$(22, 7, 360)$ & 3960 & 3 & $M_{11}$\\ \hline
$2$-$(22, 7, 720)$ & 7920 & 2 & $M_{11}$\\ \hline
$3$-$(22, 7, 18)$ & 792 & 1 & $M_{11}\times Z_2$\\ \hline
$3$-$(22, 7, 90)$ & 3960 & 3 & $M_{11}\times Z_2$\\ \hline
$3$-$(22, 7, 180)$ & 7920 & 1 & $M_{11}$\\ \hline
\hline
\end{tabular} \end{footnotesize}\caption{$t$-designs constructed from the group $M_{11}$, $v\leq 22$}\label{tb:td22}\end{center}
\end{table}

\begin{rem} 
We proved the existence of $3$-$(22,7,18)$ design, since it is the first known example of the design with these parameters. The design with parameters $5$-$(12,6,6)$ is the extension of the $4$-$(11,5,6)$ design.  The designs have the same parameters as the copies of Steiner systems $S(4,5,11)$ and $S(5,6,12)$, respectively, but the ones that we obtained are simple. Since there are no designs with parameters $5$-$(12,6,6)$ and $4$-$(11,5,6)$ mentioned in \cite{crc-t-des, crc_kreher} we conclude that the designs presented in this paper are the first examples of the simple designs with these parameters. 

$2$-designs defined on 22 points from the Table \ref{tb:td22} are not mentioned in \cite{crc-block-des} since $r>41$. Up to our best knowledge they have not been known before, so we proved their existence. The Steiner system $4$-$(11,5,1)$ is known as Witt design $W_{11}$. For further information on $W_{11}$ we refer the reader to \cite{witt1,witt2}. 
All others transitive $t$-designs described in Table \ref{tb:td22} were previously known. For further information on quasi-symmetric $3$-$(12,6,2)$ we refer the reader to \cite{quasi_symm} and for others known designs mentioned in Table \ref{tb:td22} see \cite{crc-t-des,crc-block-des}.
\end{rem}

\subsection{Block designs with $v=55$}

In this section we give all
$t$-designs with 55 points on which the group $M_{11}$ acts transitively. The
designs are obtained from the group $M_{11}$ by using Theorem \ref{main}. In that case, the stabilizer of a point is subgroup of $M_{11}$ having the index 55. 
The list of all designs obtained is given in Table \ref{tb:d55_1}. 
In Table \ref{tb:d55_1} we give the parameters of the constructed structures, the number of non-isomorphic structures and their full automorphism group.

\begin{table}[H]
\begin{center} \begin{scriptsize}
\begin{tabular}{|c|c|c|}
\hline
Parameters of block designs & \# non-isomorphic  & Full automorphism group  \\
\hline
\hline
$2$-$(55, 3, 4)$ &  1 & $M_{11}$ \\ \hline
$2$-$(55, 4, 8)$ & 1 & $M_{11}$\\ \hline
$2$-$(55, 4, 16)$ & $\geq 2$ & $M_{11}$\\ \hline
$2$-$(55, 6, 10)$ & 1 &$M_{11}$\\ \hline
$2$-$(55, 6, 20)$ & 2 & $M_{11}$\\ \hline
$2$-$(55, 6, 40)$ & $\geq 14$ & $M_{11}$\\ \hline
$2$-$(55, 7, 14)$ & 1 & $M_{11}$\\ \hline
$2$-$(55, 7, 28)$ & 5 & $M_{11}$\\ \hline
$2$-$(55, 7, 56)$ & $\geq 8$ & $M_{11}$\\ \hline
$2$-$(55, 9, 32)$ & 10 & $M_{11}$\\ \hline
$2$-$(55, 9, 48)$ & 23 & $M_{11}$\\ \hline
$2$-$(55, 9, 64)$ & 16 & $M_{11}$\\ \hline
$2$-$(55, 9, 96)$ & $\geq 45$ & $M_{11}$\\ \hline
$2$-$(55, 10, 20)$ & 1 & $M_{11}$\\ \hline
$2$-$(55, 10, 40)$ & 8 & $M_{11}$\\  \hline
$2$-$(55, 10, 48)$ & 6 & $M_{11}$\\ \hline
$2$-$(55, 10, 60)$ & 27 & $M_{11}$\\ \hline
$2$-$(55, 10, 80)$ & 9 &  $M_{11}$\\  \hline
$2$-$(55, 10, 120)$ & $\geq 268$ &  $M_{11}$\\  \hline
$2$-$(55, 12, 44)$ & 5 & $M_{11}$\\ \hline
$2$-$(55, 12, 88)$ & 43 & $M_{11}$\\ \hline
$2$-$(55, 12, 176)$ & $\geq 767$ & $M_{11}$\\ \hline
$2$-$(55, 13, 104)$ & 81 &   $M_{11}$\\ \hline
$2$-$(55, 13, 208)$ & $\geq 498$ &   $M_{11}$\\ \hline
$2$-$(55, 15, 112)$ & 8 &   $M_{11}$\\ \hline
$2$-$(55, 15, 140)$ & 53 & $M_{11}$\\ \hline
$2$-$(55, 15, 280)$ & $\geq 559$ & $M_{11}$\\ \hline
$2$-$(55, 16, 80)$ & 3 & $M_{11}$\\ \hline
$2$-$(55, 16, 160)$ & 93 & $M_{11}$\\ \hline
$2$-$(55, 16, 320)$ & $\geq 1483$ & $M_{11}$\\ \hline
$2$-$(55, 18, 68)$ & 2 & $M_{11}$\\ \hline
$2$-$(55, 18, 102)$ & 8 &   $M_{11}$\\ \hline
$2$-$(55, 18, 136)$ & 39 &  $M_{11}$\\ \hline
$2$-$(55, 18, 204)$ & 215 &  $M_{11}$\\ \hline
$2$-$(55, 18, 272)$ & 161 & $M_{11}$\\ \hline
$2$-$(55, 18, 408)$ & $\geq 4721$ & $M_{11}$\\ \hline
$2$-$(55, 19, 152)$ & 18 & $M_{11}$\\ \hline
$2$-$(55, 19, 456)$ & $\geq 5789$ & $M_{11}$\\ \hline
$2$-$(55, 19, 228)$ & 173 &$M_{11}$\\ \hline
$2$-$(55, 19, 304)$ & 376 &  $M_{11}$\\ \hline
$2$-$(55, 21, 140)$ & 9 &   $M_{11}$\\ \hline
$2$-$(55, 21, 280)$ & 249 &   $M_{11}$\\ \hline
$2$-$(55, 21, 560)$ & $\geq 4473$ &   $M_{11}$\\ \hline
$2$-$(55, 22, 308)$ & 220 & $M_{11}$\\ \hline
$2$-$(55, 22, 616)$ & $\geq 7191$ & $M_{11}$\\ \hline
$2$-$(55, 24, 184)$ & 8 & $M_{11}$\\ \hline
$2$-$(55, 24, 368)$ & 256 & $M_{11}$\\ \hline
$2$-$(55, 24, 736)$ & $\geq 2164$ & $M_{11}$\\ \hline
$2$-$(55, 25, 200)$ & 8 & $M_{11}$\\ \hline
$2$-$(55, 25, 320)$ & 25 &   $M_{11}$\\ \hline
$2$-$(55, 25, 400)$ & 445 &  $M_{11}$\\ \hline
$2$-$(55, 25, 800)$ & $\geq 5878$ &  $M_{11}$\\ \hline
$2$-$(55, 27, 78)$ & 1 &  $S_{11}$\\ \hline
$2$-$(55, 27, 234)$ & 8 & $M_{11}$\\ \hline
$2$-$(55, 27,312)$ & 56 & $M_{11}$\\ \hline
$2$-$(55, 27,468)$ & 308 &$M_{11}$\\ \hline
$2$-$(55, 27,624)$ & 626 &  $M_{11}$\\ \hline
$2$-$(55, 27,936)$ & $\geq 8393$ &  $M_{11}$\\ \hline\hline

\end{tabular} \end{scriptsize}\caption{Block designs constructed from $M_{11}$, $v=55$ \label{tb:d55_1}}\end{center} 
\end{table}

\begin{rem}
Designs from Table \ref{tb:d55_1} are not mentioned in \cite{crc-block-des} since $r>41$. Up to our best knowledge they have not been known before, so we proved the existence of $2$-designs with the parameters listed in Table \ref{tb:d55_1}. 
\end{rem}

\subsection{Block designs with $v=66$}

In this section we give all
$t$-designs with 66 points on which the group $M_{11}$ acts transitively. The
designs are obtained from the group $M_{11}$ by using Theorem \ref{main}. In that case, the stabilizer of a point is subgroup of $M_{11}$ having the index 66. 
The list of all designs obtained is given in Table \ref{tb:d66_1}. 
In Table \ref{tb:d66_1} we give the parameters of the constructed structures, the number of non-isomorphic structures and their full automorphism group.

\begin{table}[H]
\begin{center} \begin{scriptsize}
\begin{tabular}{|c|c|c|}
\hline
Parameters of block designs & \# non-isomorphic  & Full automorphism group  \\
\hline
\hline
$2$-$(66, 13, 36)$ &  1 & $M_{11}$ \\ \hline
$2$-$(66, 13, 48)$ & 13 & $M_{11}$\\ \hline
$2$-$(66, 13, 72)$ & 43 &$M_{11}$\\ \hline
$2$-$(66, 13, 96)$ & 79 & $M_{11}$\\ \hline
$2$-$(66, 13, 144)$ & $\geq 3960$ & $M_{11}$\\ \hline
$2$-$(66, 14, 56)$ & 6 & $M_{11}$\\ \hline
$2$-$(66, 14, 84)$ & 33 & $M_{11}$\\ \hline
$2$-$(66, 14, 112)$ & 105 & $M_{11}$\\ \hline
$2$-$(66, 26, 100)$ & 2 & $M_{11}$\\ \hline
$2$-$(66, 26, 120)$ & 2 & $M_{11}$\\ \hline
$2$-$(66, 26, 200)$ & 33 & $M_{11}$\\ \hline
$2$-$(66, 26, 240)$ & 4 & $M_{11}$\\  \hline
$2$-$(66, 26, 300)$ & 159 & $M_{11}$\\ \hline
$2$-$(66, 26, 400)$ & 1799 & $M_{11}$\\  \hline
$2$-$(66, 27, 36)$ & 1 & $M_{11}$\\ \hline
$2$-$(66, 27, 72)$ & 1 & $M_{11}$\\ \hline
$2$-$(66, 27, 216)$ & 60 &  $M_{11}$\\  \hline
$2$-$(66, 27, 324)$ & 49 & $M_{11}$\\ \hline
$2$-$(66, 27, 432)$ & 412 & $M_{11}$\\ \hline
\hline
\end{tabular} \end{scriptsize}\caption{Block designs constructed from $M_{11}$, $v=66$ \label{tb:d66_1}}\end{center} 
\end{table}

\begin{rem}

Designs from Table \ref{tb:d66_1} are not mentioned in \cite{crc-block-des} since $r>41$. Up to our best knowledge they have not been known before, so we proved the existence of $2$-designs with the parameters listed in Table \ref{tb:d66_1}. 
\end{rem}

\subsection{SRGs from $M_{11}$}

Using the method described in Theorem \ref{main} and Corollary \ref{main-graph}, we obtained all the regular graphs on which the alternating group $M_{11}$ acts transitively and with at most 450 vertices. 
Using the computer search we obtained strongly regular graphs on 55, 66, 144 or 330 vertices. Finally, we determined the full automorphism groups of the constructed SRGs. 

\begin{tm} \label{srg-M11}
Up to isomorphism there are exactly $5$ strongly regular graphs with at most $450$ vertices, admitting a transitive action of the group $M_{11}$.
These strongly regular graphs have parameters $(55,18,9,4)$, $(66,20,10,4)$, $(144,55,22,20)$, $(144,66,30,30)$ and $(330,63,24,9)$.
Details about the obtained strongly regular graphs are given in Table \ref{tb:srgM11}. 
\end{tm}

\begin{table}[H]
\begin{center} \begin{footnotesize}
\begin{tabular}{|c|c|c|}
\hline
Graph $\Gamma$ & Parameters  & $Aut (\Gamma) $  \\
\hline
\hline
$\Gamma_{1}=\Gamma(M_{11},H_{36})$ & (55,18,9,4) &  $S_{11}$\\ 
$\Gamma_{2}=\Gamma(M_{11},H_{34})$ & (66,20,10,4) &  $S_{12}$\\ 
$\Gamma_{3}=\Gamma(M_{11},H_{13})$ & (144,55,22,20) &  $M_{11}$\\ 
$\Gamma_{4}=\Gamma(M_{11},H_{13})$ & (144,66,30,30) &  $M_{12}:Z_2$\\ 
$\Gamma_{5}=\Gamma(M_{11},H_{24})$ & (330,63,24,9) &  $S_{11}$\\ 
\hline
\hline
\end{tabular}\end{footnotesize} 
\caption{\footnotesize SRGs constructed from the Mathieu group $M_{11}$}\label{tb:srgM11}
\end{center}
\end{table}

\begin{rem} 
The graphs $\Gamma_{1}$ and $\Gamma_{2}$ are the triangular graphs $T(11)$ and $T(12)$, respectively. Strongly regular graphs $\Gamma_{3},\Gamma_{4}, \Gamma_{5}$ were known before. For further information we refer the reader to \cite{crc-srg,aeb}.
\end{rem}

\vspace*{0.2cm}

\noindent {\bf Acknowledgement} \\
This work has been fully supported by {\rm C}roatian Science Foundation under the project 1637. 


\end{document}